\documentclass[10pt]{article}

\usepackage[utf8]{inputenc}
\usepackage{graphicx}
\usepackage{algorithm}
\usepackage{algorithmic}
\usepackage{caption}
\usepackage{subcaption}
\usepackage{siunitx} %Pour écrire avec la notation scientifique (Ex.: \num{2e+9})
\usepackage[table]{xcolor}
\usepackage{amsmath}
\usepackage{amsfonts}
\usepackage[a4paper, left=35mm,right=35mm,top=34mm,bottom=34mm]{geometry}
\usepackage[utf8]{inputenc}
\usepackage[T1]{fontenc}
\usepackage[english]{babel}

\usepackage{fancyhdr}

\graphicspath{{/}}

\title{Multilayer Perceptron-based Surrogate Models for Finite Element Analysis}

\author{Lawson Oliveira Lima$^*$,
Julien Rosenberger$^*$,
Esteban Antier$^*$,
Fr\'ed\'eric Magoul\`es\thanks{Universit\'e Paris-Saclay, CentraleSup\'elec, MICS, Gif-sur-Yvette, France, Email: frederic.magoules@hotmail.com}
}

\begin{document}

\maketitle
\thispagestyle{fancy}

\begin{abstract}
Many Partial Differential Equations (PDEs) do not have analytical solution, and can only be solved by numerical methods.
In this context, Physics-Informed Neural Networks (PINN) have become important in the last decades, since it uses a neural network and physical conditions to approximate any functions.
This paper focuses on hypertuning of a PINN, used to solve a PDE.
The behavior of the approximated solution when we change the learning rate or the activation function (sigmoid, hyperbolic tangent, GELU, ReLU and ELU) is here analyzed.
A comparative study is done to determine the best characteristics in the problem, as well as to find a learning rate that allows fast and satisfactory learning.
GELU and hyperbolic tangent activation functions exhibit better performance than other activation functions.
A suitable choice of the learning rate results in higher accuracy and faster convergence.
\end{abstract}


\section{Introduction}

Artificial intelligence draws people attention thanks to its various spectrum of application fields such as language processing, computer vision or decision making. 
Lagaris described in \cite{Lagaris-1997} a general method to solve a PDE with a neural network instead of classical methods as Finite Element Methods (FEM). % or Finite Difference Methods.
Finite Element Methods \cite{Cia2002} \cite{Hug2003} are powerfull methods to approximate solutions of a PDE.
The accuracy of these methods can be improved using various stabilization techniques like \cite{Rus2006}, \cite{HM2004WM}, \cite{MZ2018CMAME}, \cite{DF2008}.
Since FEM solution requires the solution of a large linear system, computation cost becomes prohibitive.
As an alternative, using neural network \cite{Blechschmidt-three-ways-2021} will help construct function from one already forced to satisfy the boundary conditions of the problem.
The method is broaden in \cite{Berg-unified-ANN-for-PDE-2019} where the straightforward method of solving a PDE on complex geometries without assisting the neural network to satisfy the boundary conditions is highlighted \cite{CINT-Galerkin-2015}.
In the 1990's, Physics-Informed Neural Networks were first introduced and show great expectations in terms of efficiency and versatility \cite{PINN-2017}.

The paper is organized as follows.
In Section~\ref{section:relwor}, we present concepts in Neural Networks, Partial Differential Equations and Physics Informed Neural Networks.
In Section~\ref{section:method}, we present the methodology followed to allow hypertuning of the parameters.
In Section~\ref{section:numerical}, we present numerical results and we conclude in Section~\ref{section:conc}.

\section{Neural networks for solving Partial Differential Equations}
\label{section:relwor}

\paragraph*{Neural Network}
We will focus on the study of Multilayer Perceptron (MLP), one of the first neural networks developed. % \cite{...}.
It consists in an input layer, a vector of $n$ components, and an output layer, a vector of $p$ components.
The input is subjected to a series of linear operations and nonlinear activations are applied trough, so called hidden layers, to get the output.
A MLP can be seen as a function: $\mathcal{N} : \mathbb{K}^n \rightarrow \mathbb{K}^p$.
If we note $g_i$ the (non-linear) activation function, and $W_i$ and $b_i$ the weight and bias matrix of the $i^{th}$-hidden layer.
Then, for a MLP with $k$ hidden layers we can write:
\begin{equation*}
    \begin{split}
        \forall x \in \mathbb{K}^n,
				\quad \quad
				\mathcal{N}(x) = g_k (b_k + W_k g_{k-1}(b_{k-1} \quad \quad \\
				\quad \quad \quad \quad \quad \quad \quad \quad \quad \quad + W_{k-1} g_{k-2}(\dots(b_1 + W_1 x))))
    \end{split} 
\end{equation*}
Through training epochs, the MLP passes through the training dataset coupled with gradient descent so that $\mathcal{N}$ approaches a function of interest.

Our goal is to solve a Partial Differential Equation (PDE) defined within a domain $\Omega$:
\begin{equation*}
    \begin{split}
        (E): \left \{ \begin{array}{ll} \mathcal{Q}(u, \nabla u, H u, \dots)(x) = 0 & \text{inside $\Omega$}\\
            \mathcal{R}(u, \nabla u, H u, \dots)(x) = f(x) & \text{on $\partial \Omega$}\end{array}\right.
    \end{split} 
\end{equation*}
where $\mathcal{Q}$ and $\mathcal{R}$ are linear operators, $\Omega$ is the domain considered, $\partial \Omega$ the boundary of the domain considered, and $u : \mathbb{K}^n \rightarrow \mathbb{K}^p$ the solution to $(E)$.

We want to replace $u$ with a neural network $\mathcal{N}$ and train it to get a solution of $(E)$.
However, this would mean training on both the domain and its boundary.
To avoid this case, as the PDE is always defined with boundary conditions, we reformulate the problem by replacing $u$ with a function $\Psi$ defined as: $\Psi(x) = A(x) + F(x)\mathcal{N}(x)$ where $A$ verifies the boundary conditions, and $F$ is null on the boundary.
As such, we unconstrained the problem, and the neural network can train inside the domain, with the boundary conditions handled by the above mentioned functions.

\begin{figure}[h]
    \centering
    \begin{subfigure}[b]{.43505\linewidth}
        \centering
        \includegraphics[width=1\linewidth]{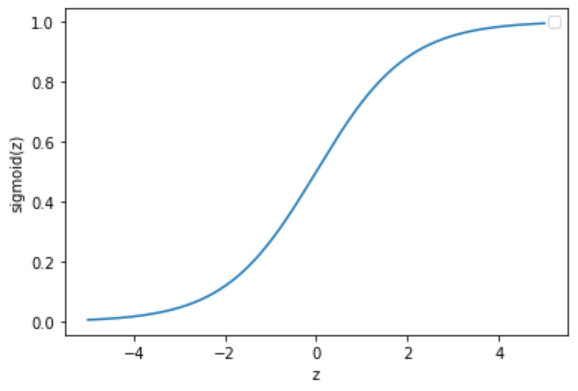}
        \caption{Sigmoid}    
    \end{subfigure}		
    \begin{subfigure}[b]{.43505\linewidth}
        \centering
        \includegraphics[width=1\linewidth]{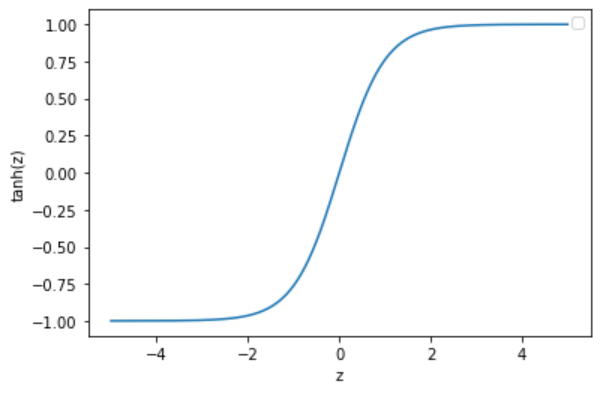}
        \caption{Hyperbolic tangent}    
    \end{subfigure}
    \begin{subfigure}[b]{.43505\linewidth}
        \centering
        \includegraphics[width=1\linewidth]{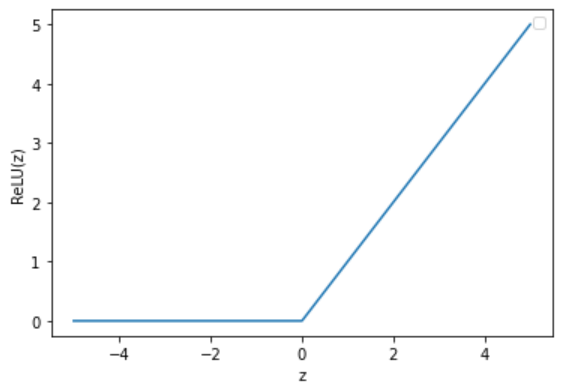}
        \caption{ReLU}
    \end{subfigure}		
    \begin{subfigure}[b]{0.43505\linewidth}
        \centering
        \includegraphics[width=1\linewidth]{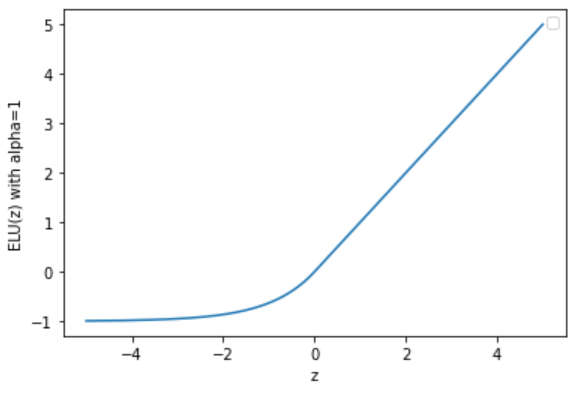}
        \caption{ELU}
    \end{subfigure}
        \begin{subfigure}[b]{0.43505\linewidth}
        \centering
        \includegraphics[width=1\linewidth]{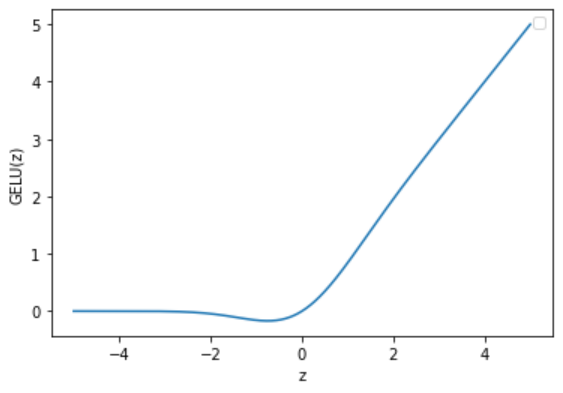}
        \caption{GELU}  
    \end{subfigure}
	\caption{Activation functions.}
\label{Activation_functions}
\end{figure}

\paragraph*{Activation functions}
%\subsection{Activation functions}
Activation functions set the behaviour of a neuron at a low level in the network.
Their impacts can be more or less fast to compute and drastically change the convergence speed and quality.
We consider in Equations~(\ref{activ_sigmoid})-(\ref{activ_gelu}) the most rampant activation functions in the state-of-the-art articles:
\begin{subequations}
    \begin{equation}
        \forall z \in \mathbb{R}, \ \sigma(z):=\frac{1}{1-e^{-z}}
        \label{activ_sigmoid}
    \end{equation}    
    \begin{equation}
        \forall z \in \mathbb{R}, \ \text{tanh}(z):=\frac{e^z-e^{-z}}{e^z+e^{-z}} = \sigma(2z)-\sigma(-2z)
        \label{activ_tanh}    
    \end{equation}
    \begin{equation}
        \alpha \in [0,1], \ \forall z \in \mathbb{R}, \ \text{ReLU}(z):= \text{max}(\alpha z,z)
        \label{activ_relu}
    \end{equation}
    \begin{equation}
        \begin{split}
        \alpha \in \mathbb{R}^*, \ \forall z \in \mathbb{R}, \ \text{ELU}(z):= \left \{ \begin{array}{ll} z , \text{if} \ z \geq 0 \\
            \alpha(e^z-1), \text{if} \ z < 0
            \end{array}\right.
        \end{split}
        \label{activ_elu}
    \end{equation} 
    \small{\begin{equation}
            \forall z \in \mathbb{R},  
             ~\text{GELU}(z)=0.5z\displaystyle\left(1+\text{tanh}\left[\sqrt{\frac{2}{\pi}}(z+0.044715z^3)\right]\displaystyle\right)
        \label{activ_gelu} 
    \end{equation}}
%\label{activ_funcs}
\end{subequations}
The first one Equation~(\ref{activ_sigmoid}) coming from an imitation of biological neurons' behaviour is the sigmoid.
Inspired from the previous one, but zero-centered, is the hyperbolic tangent Equation~(\ref{activ_tanh}).
Yet, activation functions easier and faster to compute arose.
The most notorious one is the Rectified Linear Unit (ReLU) function Equation~(\ref{activ_relu}) and it holds many variants. % such as LeakyReLU.
A smoother variant is the Exponential Linear Unit (ELU) function.
This function is defined Equation~(\ref{activ_elu}).
The final variant to introduce is the Gaussian Error Linear Unit (GELU) function presented in \cite{GELU-2018}.
This function tries to embed the Dropout regularizer by randomly multiplying outputs of neurons to the activation function itself.
For using this function, the inputs have to be batch normalized first.
The GELU function is defined Equation~(\ref{activ_gelu}).
The different activation functions are represented in Figure~\ref{Activation_functions}.

%\vspace{-0.5cm}

\section{Methodology}
\label{section:method}

\paragraph*{Problem Statement}
We want to solve the following PDE in a square domain $\Omega = {[0,1]}^2$ using Physics-Informed Neural Networks (PINN).
\begin{equation}
    \Delta \psi(x,y) + \psi(x,y)\cdot\frac{\partial \psi(x,y)}{\partial y} = f(x,y) \quad\text{on}\quad \Omega = {[0,1]}^2  
    \label{pde_example}
\end{equation}
where $f(x, y)=\sin(\pi x)(2-\pi^2y^2+2y^3\sin(\pi x))$
and with the boundary conditions:
\begin{equation*}
    \begin{split}
    \left\{ \begin{array}{lll} \psi(0,y) = \psi(1,y) = \psi(x,0)=0 \\
        \frac{\partial \psi}{\partial y}(x,1) = 2\sin(\pi x) 
        \end{array}\right.
    \end{split}
\end{equation*}   
This problem has an analytical solution:
$\psi_{th}(x, y) = y^2\sin(\pi x).$
Thus, to solve it, we define our loss function as:
\begin{equation*}
    \mathcal{L} = ||\Delta \psi(x,y) + \psi(x,y)\cdot\frac{\partial \psi(x,y)}{\partial y} - f(x,y) ||_2 
\end{equation*} 
with $\psi_{ap}(x,y)$ defined as
    $\psi_{ap}(x,y) = A(x,y) + F(x,y) N(x,y)$
and with $A$ and $F$ equal to
    $A(x,y) = y\sin(\pi x)$
   and $F(x,y) = \sin(x-1)\sin(y-1)\sin(x)\sin(y)$.
A careful reader will notice the lack of knowledge about the exact solution of the PDE.
Furthermore, the model is called a PINN as the loss takes into account the residual of Equation~(2).

\paragraph*{Implementation}

The implementation of the multilayer perceptron is similar to \cite{MZE2013EB} \cite{MZ2012ER} and is done using the Python library Jax, 
which has high performance and low computational cost due to being built using XLA, and is thus recommended for simulations involving neural networks.
Furthermore, it can be used to vectorize functions and execute them faster in Graphics Processing Unit (GPU).
The steps of the model used are described in the algorithm below.
\begin{algorithm}[h]
    \caption*{Solving PDEs using neural networks (PINN method)}
    \hspace*{\algorithmicindent} \textbf{Input} domain and boundary-coords, tolerance, $A$ and $F$
    \hspace*{\algorithmicindent} \textbf{Output} trained-PINN
    \begin{algorithmic}
    	  \STATE set max-epoch, max-train-step, batch-size, validation-size
    	  \STATE set learning rate and activation function		
        \WHILE {loss $\geq$ tolerance}
    	    \FOR{epoch $\leq$ max-epoch}
    	        \STATE sample-coords $\leftarrow$ sample batch-size points + noise
    	        \STATE loss  $\leftarrow$ compute $\mathcal{L}2$-norm on PINN-model and the PDE
    	        \STATE apply gradient descent and backpropagation to reduce loss
    	    \ENDFOR
        \ENDWHILE
        \STATE make the validation of the model
    	\RETURN PINNmodel
    \end{algorithmic}
\end{algorithm}

To achieve a better performance in a shorter time, we have vectorized it to compute in parallel the neural network evaluation and make gradient descent on GPU.

\paragraph*{Hyperparameters}
The hypertuning comes when one wants to construct a model and a training that will end to a low validation error.
Let's list all the hyperparameters involved in this scheme:
(i) \texttt{l\_units}: list of integers depicting the number of hidden layers and the number of neurons per layer of the sequential machine learning model;
(ii) \texttt{l\_activations}: list depicting the activation functions of each layer of the model;
(iii) \texttt{noise}: float in $[0,1]$ conditioning the use of a Gaussian layer of mean $0$ after the input layer;
(iv) \texttt{stddev}: the standard deviation of the Gaussian layer when it is used;
(v) \texttt{optimizer}: string depicting the optimizer used;
(vi) \texttt{learning\_rate}: float setting the learning rate of the previous optimizer;
(vii) \texttt{epochs\_max}: integer setting the number of maximum epoch loops;
(viii) \texttt{batch\_size}: integer setting the number of points extracted from the mesh to construct the training set at each epoch loop.
Hypertuning consists of adapting the parameters of the neural network in order to obtain a model that has the lowest validation error.
To do this, several models are trained and the one with the best performance is chosen

\section{Numerical experiments}
\label{section:numerical}

For the hyperparameters, we used $batch\_size = 50$ and a multilayer perceptron with 3 completely connected layers, the first layer has 2 neurons, the second 30 and the last 1, each followed by the activation function.
Thus the neural network considered has 121 trainable parameters: hidden layer = $2 * 30+30=90$ and output layer = $30* 1 + 1=31$.

\subsection{First parameter: Learning rate}

\begin{figure}[h]
    \centering
    \begin{subfigure}[b]{0.43505\linewidth}
    \centering
    \includegraphics[width=1\linewidth]{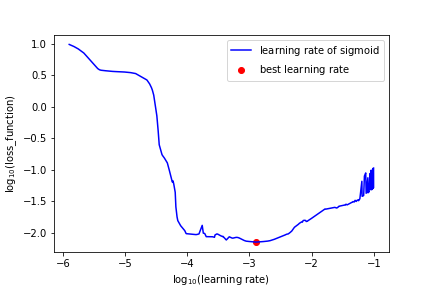}
    \caption{Sigmoid}
    \end{subfigure}
    \begin{subfigure}[b]{0.39505\linewidth}
    \centering
    \includegraphics[width=1\linewidth]{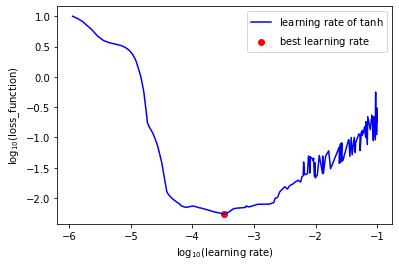}
    \caption{Tanh}
    \end{subfigure}
    \begin{subfigure}[b]{0.43505\linewidth}
    \centering
    \includegraphics[width=1\linewidth]{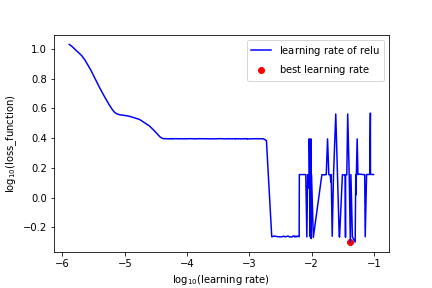}
    \caption{ReLU}
    \end{subfigure}		
    \begin{subfigure}[b]{0.43505\linewidth}
    \centering
    \includegraphics[width=1\linewidth]{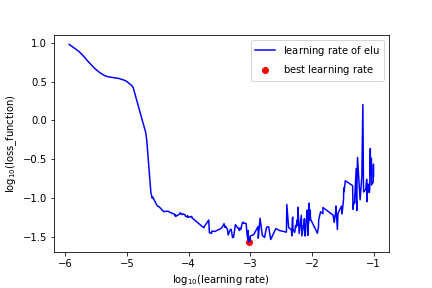}
    \caption{ELU}
    \end{subfigure}
    \begin{subfigure}[b]{0.43505\linewidth}
    \centering
    \includegraphics[width=1\linewidth]{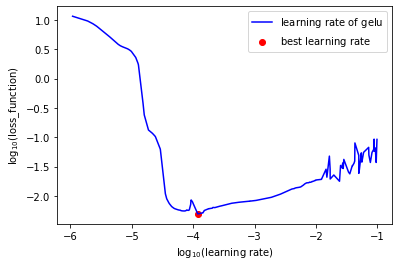}
    \caption{GELU}
    \end{subfigure}
    \caption{Learning rates.}
    \label{spaces}
    \end{figure}

Studying Equation~(\ref{pde_example}) we note that each activation function has a learning rate that improve its accuracy and reduce the loss function.
Thus, to get a good solution, we have to find a model that minimizes it.
It's possible to see the behavior between learning rate and loss function in Figure~\ref{spaces} for each one of the studied functions.
In order to find it, we have used grid and random search in pre-defined intervals.

\paragraph*{Grid search}
Grid search is a very common way to make hyperparameter tuning.
It consists of determining a domain and then using points that are uniformly distributed.
So, to use this method, we chose 5 intervals with 50 learning rates in each and then ran simulations to obtain the best loss function.
The mean absolute error and respective learning rate for each activation function using grid search  are presented in Table~\ref{grid_search}.
\begin{table}[h]
    \centering
    \resizebox{9cm}{!}{
    \small
    \begin{tabular}{cccccc}
        \hline
        MAE at epoch & Tanh & Sigmoid & GELU & ELU & ReLU \\
        \hline %\hline
        $10 000$ &  \cellcolor[HTML]{d3d3d3} $1.71\cdot 10^{-5}$ & $3.10\cdot 10^{-5}$ &  $1.84\cdot 10^{-4}$ & $5.38\cdot 10^{-5}$ & $0.016$ \\
        \hline
        $20 000$ & $1.43\cdot 10^{-5}$ & $2.18\cdot 10^{-5}$ & \cellcolor[HTML]{d3d3d3} $9.92\cdot 10^{-6}$ & $1.65\cdot 10^{-4}$ & $0.016$ \\
        \hline
        $50 000$ &  $1.05\cdot 10^{-5}$  & $1.08\cdot 10^{-5}$ & \cellcolor[HTML]{d3d3d3} $8.03\cdot 10^{-6}$ & $1.57\cdot 10^{-4}$ & $0.016$ \\
        \hline
        Learning rate &  $3.2040\cdot10^{-4}$ & $1.3673\cdot10^{-5}$ & $1.1837\cdot10^{-4}$ & $1.1837\cdot10^{-3}$ & $6.6939\cdot10^{-2}$ \\
        \hline
    \end{tabular}}
    \caption{Grid search (grey is the best accuracy).}
\label{grid_search}
\end{table}

\paragraph*{Random search} 
Random search is similar in principle to grid search, using points distributed over a domain.
However, the main difference is that they are spread randomly, which allows testing parameters that would not be tested using the first method.
Thus, this procedure was performed on the same intervals as previously used.
Results of the mean absolute error and respective best learning rate for each activation function using random search are collected in Table~\ref{random_search}.
\begin{table}[h]
    \centering
    \resizebox{9cm}{!}{
    \small
    \begin{tabular}{cccccc}
        \hline
        MAE at epoch & Tanh & Sigmoid & GELU & ELU & ReLU \\
        \hline 
        $10 000$ &  \cellcolor[HTML]{d3d3d3} $3.10\cdot 10^{-5}$ & $6.78\cdot 10^{-5}$ &  $1.85\cdot 10^{-4}$ & $5.48\cdot 10^{-5}$ & $0.016$ \\
        \hline
        $20 000$ & \cellcolor[HTML]{d3d3d3} $2.18\cdot 10^{-5}$ & $5.75\cdot 10^{-5}$ &  $1.47\cdot 10^{-5}$ & $5.86\cdot 10^{-5}$ & $0.016$ \\
        \hline
        $50 000$ &  $1.12\cdot 10^{-5}$  & $1.08\cdot 10^{-5}$ & \cellcolor[HTML]{d3d3d3} $6.31\cdot 10^{-6}$ & $1.08\cdot 10^{-4}$ & $0.016$ \\
        \hline
        Learning rate &  $3.2915\cdot10^{-4}$ & $1.2626\cdot10^{-5}$ & $1.1891\cdot10^{-4}$ & $9.5295\cdot10^{-4}$ & $4.1886\cdot10^{-2}$ \\
        \hline
    \end{tabular}}
    \caption{Random search (grey is the best accuracy).}
\label{random_search}
\end{table}

\paragraph*{Analysis} 
Mean absolute error and respective learning rate for each activation function with best learning rates are presented in Table~\ref{best_search}.
\begin{table}[h]
    \centering
    \resizebox{9cm}{!}{
    \small
    \begin{tabular}{cccccc}
        \hline
        MAE at epoch & Tanh & Sigmoid & GELU & ELU & ReLU \\
        \hline 
        $10 000$ &  \cellcolor[HTML]{d3d3d3} $1.71\cdot 10^{-5}$ & $3.10\cdot 10^{-5}$ &  $1.85\cdot 10^{-4}$ & $5.48\cdot 10^{-5}$ & $0.016$ \\
        \hline
        $20 000$ & \cellcolor[HTML]{d3d3d3} $1.43\cdot 10^{-5}$ & $2.18\cdot 10^{-5}$ &  $1.47\cdot 10^{-5}$ & $5.86\cdot 10^{-5}$ & $0.016$ \\
        \hline
        $50 000$ &  $1.12\cdot 10^{-5}$  & $1.05\cdot 10^{-5}$ & \cellcolor[HTML]{d3d3d3} $6.31\cdot 10^{-6}$ & $1.08\cdot 10^{-4}$ & $0.016$ \\
        \hline
        Learning rate &  $3.2040\cdot10^{-4}$ & $1.3673\cdot10^{-5}$ & $1.1891\cdot10^{-4}$ & $9.5295\cdot10^{-4}$ & $4.1886\cdot10^{-2}$ \\
        \hline
    \end{tabular}}
    \caption{Best learning rates (grey is the best accuracy).}
\label{best_search}
\end{table}
We can see that random search presented the lowest errors for almost all functions, except for sigmoid.
An expected result, since it does not follow a pattern to test each learning rate.
Furthermore, it can be seen that the average absolute error obtained is quite small for GELU, tanh and sigmoid. 
We compare for the same accuracy with FEM on GPU.
The linear system of equations is solved with the Alinea library \cite{MC2015IJHPCA}.
The problem is reformulated in parallel with a domain decomposition method \cite{MCP2016IJCM} and is solved on GPU \cite{MC2016CC} \cite{CM2017JOS}.
The local solver inside each subdomain is performed with the conjugate gradient method \cite{CM2017AES} involving a auto-tuning of the GPU memory \cite{MCP2015}.
We obtain a computational time of 25 seconds for $50 \ 10^3$ epochs for the MLP surrogate models and for the FEM method involving $2.5 \ 10^6$ nodes, more time.

\subsection{Second parameter: Activation function}

\begin{figure}[h]
    \centering
    \begin{subfigure}[b]{0.43505\linewidth}
        \centering
        \includegraphics[width=1\linewidth]{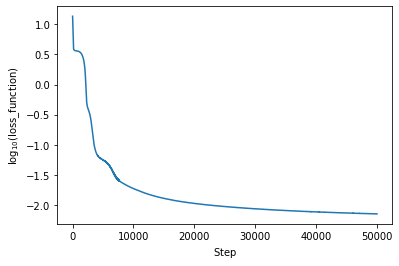}
        \caption{Training (sigmoid)}
    \end{subfigure}
    \begin{subfigure}[b]{.43505\linewidth}
        \centering
        \includegraphics[width=1\linewidth]{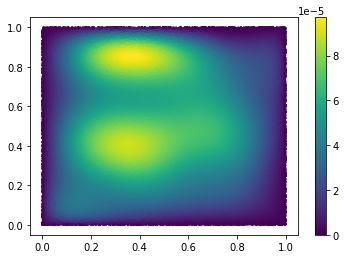}
        \caption{Absolute error (sigmoid)}    
    \end{subfigure}
    \begin{subfigure}[b]{0.43505\linewidth}
        \centering
        \includegraphics[width=1\linewidth]{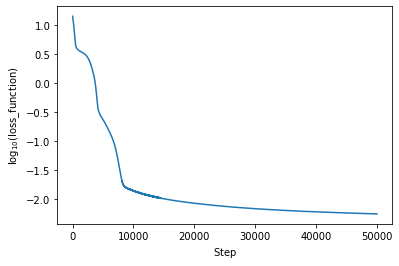}
        \caption{Training (hyperbolic tangent)}
    \end{subfigure}
    \begin{subfigure}[b]{.43505\linewidth}
        \centering
        \includegraphics[width=1\linewidth]{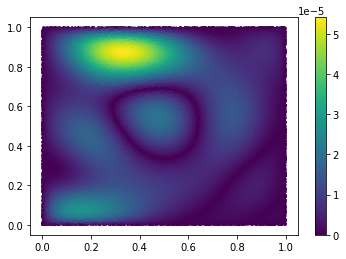}
        \caption{Absolute error (hyperbolic tangent)}    
    \end{subfigure}
    \begin{subfigure}[b]{0.43505\linewidth}
        \centering
        \includegraphics[width=1\linewidth]{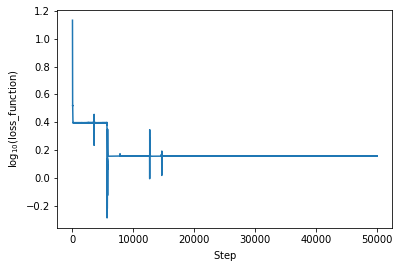}
        \caption{Training (ReLU)}
    \end{subfigure}
    \begin{subfigure}[b]{.43505\linewidth}
        \centering
        \includegraphics[width=1\linewidth]{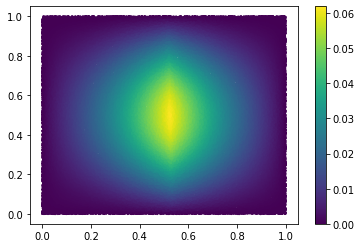}
        \caption{Absolute error (ReLU)}    
    \end{subfigure}
    \begin{subfigure}[b]{0.43505\linewidth}
        \centering
        \includegraphics[width=1\linewidth]{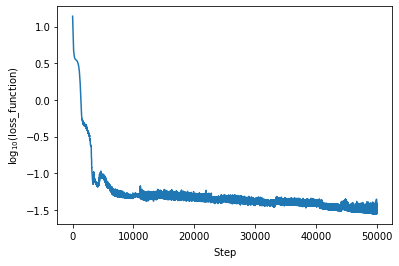}
        \caption{Training (ELU)}
    \end{subfigure}
    \begin{subfigure}[b]{.43505\linewidth}
        \centering
        \includegraphics[width=1\linewidth]{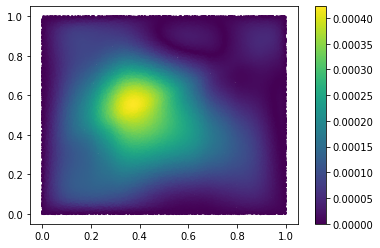}
        \caption{Absolute error (ELU)}    
    \end{subfigure}
    \begin{subfigure}[b]{0.43505\linewidth}
        \centering
        \includegraphics[width=1\linewidth]{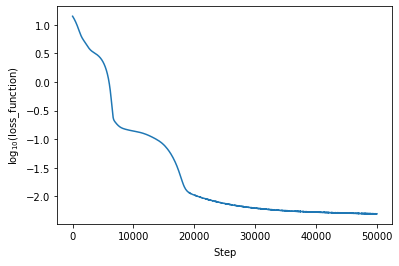}
        \caption{Training (GELU)}
    \end{subfigure}
    \begin{subfigure}[b]{.43505\linewidth}
        \centering
        \includegraphics[width=1\linewidth]{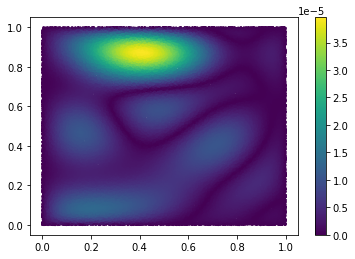}
        \caption{Absolute error (GELU)}    
    \end{subfigure}
    \caption{Training and absolute error.}
\label{PDE8}
\end{figure}

With the best learning rate for each activation function, we report in Figure~\ref{PDE8}, the comparison in order to find the best solution to the problem.
The activation functions ELU and ReLU did not obtain satisfactory results since they reached a local minimun that did not allow for good accuracy.
On the other hand, GELU has the best performance and the hyperbolic tangent has similar results.
Finally, we see that sigmoid has performance as good as GELU and tanh.
Additional simulations with varying hyperparameters did not change the performance of one over the other.

\section{Conclusions}
\label{section:conc}

It was observed that GELU and hyperbolic tangent has lower absolute error for a smaller number of epochs, a very favorable characteristic, since it makes the neural network easier to train.
Also, sigmoid performed well, but ELU and ReLU did not obtain satisfactory results.
Moreover, with the optimization of the learning rate in the hyperbolic tangent and GELU a good mean absolute error was obtained for 10 000 and 20 000 epochs, which enable us to reduce the number of epochs and the training time.

\section{Acknowledgements}

The authors would like to thank
M.~Besbes, A.~Paun, and G.~Ruault for the useful discussions.

\end{document}